\documentclass[preprint,12pt]{elsarticle}

\usepackage{amsmath,amsthm,amssymb}

\usepackage{bm}
\usepackage{mathrsfs}
\usepackage{enumitem}
\usepackage{graphicx}
\usepackage{bbm}
\usepackage{color}
\usepackage{geometry}
\newtheorem{theorem}{Theorem}
\newtheorem{lemma}{Lemma}

\newtheorem{definition}{Definition}

\numberwithin{equation}{section}

\journal{}

\begin{document}

\begin{frontmatter}



\title{Diffusion approximations for interacting stochastic systems with reflection and control}


\author[mymainaddress]{Thoa Thieu \corref{mycorrespondingauthor}}
\cortext[mycorrespondingauthor]{Corresponding author}
\author[mysecondaryaddress]{Roderick Melnik}

\address[mymainaddress]{School of Mathematical and Statistical Sciences, The University of Texas Rio Grande Valley,\\ 1201 W University Dr, Edinburg, TX 78539, USA}
\address[mysecondaryaddress]{MS2Discovery Interdisciplinary Research Institute, Wilfrid Laurier University, \\75 University Ave W, Waterloo, Ontario, Canada N2L 3C5 \\Emails: thoa.thieu@utrgv.edu, rmelnik@wlu.ca}


\begin{abstract}
We study diffusion approximations for a class of interacting stochastic systems with reflection and control. Motivated by interacting stochastic dynamics subject to feedback mechanisms and boundary constraints, we consider diffusion-scaled stochastic processes incorporating stochastic fluctuations, state-dependent interactions, and reflection. Under suitable assumptions, we establish convergence in distribution of the scaled processes to systems of interacting reflected stochastic differential equations of Ornstein-Uhlenbeck type. The limiting dynamics capture key features of constrained multi-agent systems, including mean-reverting behavior, interaction effects, and confinement within bounded domains through Skorokhod reflection. The analysis combines diffusion-scaling arguments, stability estimates, and continuity properties of the Skorokhod map to connect discrete stochastic systems with their reflected diffusion limits. To illustrate the framework, we present numerical examples motivated by crowd dynamics and neural population dynamics. The simulations demonstrate qualitative agreement between the finite stochastic systems and the corresponding reflected diffusion models and illustrate how diffusion approximations can provide tractable descriptions of interacting stochastic systems with constraints.

\end{abstract}


%
%
%

%
%
%

\end{frontmatter}



	\section{Introduction}
	
	Modeling the movement and interaction of multiple agents in constrained environments is a fundamental problem in a wide range of applications, including operations research, crowd dynamics, traffic systems, and biological processes. In many such settings, agents evolve under a combination of stochastic fluctuations, individual objectives, and mutual interactions, while remaining confined to a bounded or structured domain. Developing mathematical models that capture these features in a tractable and flexible manner is therefore of considerable interest.
	
	In the context of crowd dynamics, individuals typically exhibit several characteristic behaviors. They tend to move toward their intended destinations, such as exits in evacuation scenarios, while simultaneously avoiding collisions with others through repulsive interactions. Their motion also takes place within spatially constrained environments, such as rooms, corridors, or open areas with obstacles. Stochastic differential equations (SDEs) provide a natural modeling framework for these phenomena. In particular, Ornstein-Uhlenbeck (OU) type processes capture mean-reverting behavior toward target locations, while additive noise models inherent variability. When multiple agents are considered, interaction effects can be incorporated through additional drift terms, leading to systems of interacting stochastic processes.
	
	In many realistic situations, the dynamics are further subject to constraints such as physical boundaries or capacity limitations. These constraints can be modeled through reflected stochastic differential equations, in which a reflection term ensures that the process remains within an admissible domain. The theory of reflected SDEs is well established. Early work by Skorokhod \cite{Skorohod1961} introduced reflected diffusion processes in one dimension, and subsequent extensions to higher dimensions and more general domains were developed in \cite{Lions1984,Saisho1987,Dupuis1993}. Further results on existence, uniqueness, and stability of reflected SDEs can be found in \cite{Slominski1993,Slominski1994,MarinRubio2004,Wells2006,Sabelfeld2019,Li2020novel,Bhattacharya2021random,Falkowski2025,Bass2025uniqueness}.
	
	{\color{black}{
			Beyond the classical theory of reflected stochastic differential equations,
			substantial progress has also been made in the study of interacting particle
			systems and mean-field models with reflection. In particular, mean-field and
			McKean-Vlasov stochastic systems have been widely used to describe collective
			dynamics in large populations of interacting agents, where the number of
			particles tends to infinity and interactions are represented through the
			empirical distribution of the system. These models have attracted considerable
			attention in recent years due to their applications in biology, economics,
			engineering, and social dynamics \cite{Carmona2013control,Liu2021long,Jiang2026learning}.
			
			The present work takes a different perspective. Rather than studying a
			mean-field limit as the number of agents tends to infinity, we focus on
			diffusion approximations arising from the scaling of underlying stochastic
			systems with reflection and control. Our objective is to derive finite-dimensional
			systems of interacting reflected diffusions from appropriately scaled stochastic
			dynamics and to establish rigorous convergence results under general assumptions.
			
	}}
	Beyond their direct formulation, reflected diffusion models can also arise as limits of appropriately scaled discrete stochastic systems. When the underlying dynamics are driven by random events occurring at high rates, it is natural to consider diffusion scaling, under which centered and normalized processes converge in distribution to continuous stochastic dynamics. This perspective provides a connection between discrete stochastic models and continuous SDE-based descriptions, and allows one to interpret reflected diffusions as tractable approximations of complex interacting systems. While such constructions are often inspired by queueing-type or birth-death processes, the resulting framework is not restricted to classical queueing networks and can accommodate a broader class of interacting stochastic systems.
	
	Queueing-type models nevertheless provide useful intuition and have been widely applied in areas such as call centers, manufacturing systems, communication networks, and traffic flow; see, for example, \cite{Bramson2006stability,Lee2011,Pilipenko2014introduction,Kumaran2019queuing,Ata2024singular,Xu2024pigat,Su2025improved}. In traffic modeling, queue-based approaches have been used to study flow dynamics and congestion effects \cite{Vandaele2000queueing,Helbing2003section}. More recent developments incorporate stochastic variability and dependence structures to better capture complex system behavior \cite{Bogachev2024approximate}. In this work, such queue-inspired constructions serve primarily as a motivation for discrete stochastic approximations rather than as a direct modeling framework.
	
	The objective of this paper is to develop and analyze a class of diffusion approximations for interacting stochastic systems with reflection and control. Starting from discrete stochastic processes that incorporate stochastic fluctuations, interaction effects, and boundary constraints, we establish conditions under which suitably scaled processes converge in distribution to systems of reflected stochastic differential equations of Ornstein--Uhlenbeck type. The analysis combines a decomposition of the prelimit processes, functional limit theorems for the driving noise, a priori estimates, and continuity properties of the Skorokhod map.
	
	The resulting framework provides a flexible approach for modeling multi-agent systems in constrained environments. It captures key features such as goal-directed motion, interaction among agents, and confinement within admissible domains, while remaining amenable to analysis through diffusion approximations.
	
	We illustrate the applicability of the approach through two numerical examples. The first considers a crowd dynamics scenario in which interacting agents move within a confined domain, while the second examines a neural population model describing stochastic activity under competition and boundary constraints. These examples demonstrate that discrete stochastic systems can exhibit behavior consistent with the corresponding reflected diffusion models under appropriate scaling, and highlight the usefulness of diffusion approximations in describing complex multi-agent dynamics.
	
	The paper is organized as follows. In Section~2, we introduce the stochastic model and formulate the reflected diffusion framework. Section~3 presents the main diffusion approximation result and its proof. Section~4 provides numerical experiments illustrating the behavior of the discrete and diffusion models. Finally, concluding remarks and future directions are discussed in Section~5.

	\section{Setting of the model equations}
	
	We consider a system of $N$ interacting stochastic processes
	${Q_t^{(n), i}}_{i=1}^N$, where each component
	\[
	Q_t^{(n), i} \in \mathbb{R}^2
	\]
	represents the state of agent $i$ in a two-dimensional domain.
	
	The dynamics are given by
	\[
	\mathbf{Q}_t^{(n)}
	= \mathbf{Q}_0^{(n)}
	+
	\mathbf{A}_t^{(n)} - 
	\mathbf{D}_t^{(n)}
	+
	\mathbf{C}_t^{(n)}
	+
	\mathbf{R}_t^{(n)},
	\qquad t \ge 0,
	\]
	where $\mathbf{Q}_t^{(n)} = (Q_t^{(n),1}, \ldots, Q_t^{(n),N}) \in \mathbb{R}^{2N}$, and:
	
	\begin{itemize}
		\item $\mathbf{Q}_0^{(n)}$ is the initial condition,
		
		\item $\mathbf{A}_t^{(n)}$ and $\mathbf{D}_t^{(n)}$ are stochastic input and output processes with intensities of order $n$, representing high-frequency random fluctuations,
		
		\item $\mathbf{C}_t^{(n)}$ is a control term that models mean-reverting behavior toward a nominal target $\boldsymbol{\mu} = (\mu_1,\ldots,\mu_N) \in \mathbb{R}^{2N}$,
		
		\item $\mathbf{R}_t^{(n)}$ is a reflection term ensuring that the process remains in a convex domain
		\[
		D^N \subset \mathbb{R}^{2N}.
		\]
	\end{itemize}
	
	We introduce the diffusion-scaled processes
	\[
	X_t^{(n), i} := \frac{Q_t^{(n), i} - n\mu_i}{\sqrt{n}},
	\qquad
	\mathbf{X}_t^{(n)} := (X_t^{(n),1}, \ldots, X_t^{(n),N}) \in \mathbb{R}^{2N}.
	\]
	
	The domain $D^N$ is assumed to be convex and may represent spatial constraints for each agent, as well as interaction constraints such as collision avoidance. The reflection term $\mathbf{R}_t^{(n)}$ is defined through a Skorokhod-type mechanism ensuring that $\mathbf{Q}_t^{(n)} \in D^N$ for all $t \ge 0$.
	{\color{black}{
			\subsection{A representative prelimit construction}
			
			To illustrate assumptions (H1)-(H5), we describe a representative class
			of interacting stochastic systems that fits within the abstract framework
			considered in this paper. The main theorem is formulated under assumptions
			(H1)-(H5) and therefore applies to a broader class of models beyond the
			specific example presented below.
			
			For each agent $i=1,\ldots,N$, let
			\[
			Q_t^{(n),i}\in \mathbb{R}^2
			\]
			denote the unscaled state of the $i$-th component. The dynamics are driven
			by high-frequency random fluctuations, a deterministic feedback control,
			and a reflection mechanism. More precisely, we write
			\[
			\mathbf{Q}_t^{(n)}
			=
			\mathbf{Q}_0^{(n)}
			+
			\mathbf{A}_t^{(n)}
			-
			\mathbf{D}_t^{(n)}
			+
			\mathbf{C}_t^{(n)}
			+
			\mathbf{R}_t^{(n)}.
			\]
			
			Here $\mathbf{A}^{(n)}$ and $\mathbf{D}^{(n)}$ are input and output counting
			processes with intensities of order $n$. For example, componentwise, one may take
			\[
			A_{\ell}^{(n)}(t)
			=
			N_{\ell,+}(n\lambda_{\ell}t),
			\qquad
			D_{\ell}^{(n)}(t)
			=
			N_{\ell,-}(n\mu_{\ell}t),
			\qquad
			\ell=1,\ldots,2N,
			\]
			where $N_{\ell,+}$ and $N_{\ell,-}$ are independent unit-rate Poisson
			processes. Then
			\[
			\frac{1}{\sqrt n}
			\left(
			\mathbf{A}^{(n)}_t
			-
			\mathbf{D}^{(n)}_t
			-
			\mathbb{E}
			\bigl[
			\mathbf{A}^{(n)}_t
			-
			\mathbf{D}^{(n)}_t
			\bigr]
			\right)
			\Rightarrow
			\Sigma \mathbf{W}_t
			\]
			by the functional central limit theorem.
			
			The control term is chosen so that, after centering and diffusion scaling,
			it produces the desired interaction drift. Namely, if
			\[
			\mathbf{X}_t^{(n)}
			=
			\frac{\mathbf{Q}_t^{(n)}-n\boldsymbol{\mu}}{\sqrt n},
			\]
			we assume that
			\[
			\frac{1}{\sqrt n}\mathbf{C}_t^{(n)}
			=
			\int_0^t
			\mathbf{b}\bigl(\mathbf{X}_s^{(n)}\bigr)\,ds
			+
			\boldsymbol{\varepsilon}_t^{(n)},
			\]
			where
			\[
			\sup_{0\le t\le T}
			\bigl|
			\boldsymbol{\varepsilon}_t^{(n)}
			\bigr|
			\longrightarrow 0
			\]
			in probability and
			\[
			b_i(x)
			=
			\theta_i(\mu_i-x_i)
			+
			\sum_{j\neq i}
			F_{ij}(x_i,x_j).
			\]
			
			Thus the control term generates the mean-reverting force toward the target
			position together with the interaction forces among the agents.
			
			The reflection term $\mathbf{R}^{(n)}$ is defined through the Skorokhod map
			on the admissible domain $D^N$. Equivalently,
			\[
			\mathbf{X}^{(n)}
			=
			\Gamma
			\left(
			\mathbf{X}_0^{(n)}
			+
			\mathbf{M}^{(n)}
			+
			\frac{1}{\sqrt n}
			\mathbb{E}
			\bigl[
			\mathbf{A}^{(n)}
			-
			\mathbf{D}^{(n)}
			\bigr]
			+
			\frac{1}{\sqrt n}
			\mathbf{C}^{(n)}
			\right).
			\]
			
			This example illustrates how the limiting reflected diffusion arises from
			the combined effects of high-frequency stochastic fluctuations, feedback
			control, interaction forces, and boundary reflection. The diffusion
			approximation theorem developed below is stated in the abstract framework
			(H1)-(H5), allowing it to apply to a broader class of interacting
			stochastic systems beyond this representative construction.

			\paragraph{Relation between the prelimit and limiting models.}
			The diffusion approximation theorem is formulated under assumptions
			(H1)-(H5). To make the underlying stochastic dynamics more transparent,
			we first present a representative prelimit construction illustrating how
			the primitive noise, control, and reflection mechanisms generate the
			limiting reflected diffusion.	}}
	
	\subsection{Preliminaries and assumptions}
	We rely on the following assumptions
	\begin{itemize}
		\item[(H1)] \textbf{Initial data.} 
		\[
		\mathbf{X}^{(n)}_0 \Rightarrow \mathbf{X}_0
		\quad\text{in } \mathbb{R}^{2N},
		\]
		for some $D^N$-valued random variable $\mathbf{X}_0$.

		\item[(H2)] \textbf{Primitive noise convergence.}
		See, for example, Lemma~\ref{lem:fclt}. Define
		\[
		\mathbf{M}^{(n)}_t
		:=
		\frac{1}{\sqrt n}
		\Big(
		\mathbf{A}^{(n)}_t-\mathbf{D}^{(n)}_t
		-\mathbb{E}[\mathbf{A}^{(n)}_t-\mathbf{D}^{(n)}_t]
		\Big).
		\]
		There exists a $2N$-dimensional Brownian motion $\mathbf{W}$ and a matrix
		$\Sigma \in \mathbb{R}^{2N\times 2N}$ such that
		\[
		\mathbf{M}^{(n)} \Rightarrow \Sigma \mathbf{W}
		\quad\text{in } D([0,T],\mathbb{R}^{2N}).
		\]
		
		\item[(H3)] \textbf{Critical-load centering.}
		For each $i=1,\dots,N$,
		\[
		\sup_{0\le t\le T}
		\left|
		\frac{1}{\sqrt n}\,
		\mathbb{E}\!\left[A_t^{(n),i}-D_t^{(n),i}\right]
		\right|
		\longrightarrow 0.
		\]
		In particular, this condition is satisfied if the primitive processes are
		asymptotically critically loaded, i.e.,
		\[
		\mathbb{E}[A_t^{(n),i}-D_t^{(n),i}]
		=
		n(\lambda_i^{(n)}-\mu_i^{(n)})t
		\quad\text{with}\quad
		\lambda_i^{(n)}-\mu_i^{(n)}=O(n^{-1/2}).
		\]
		
		\item[(H4)] \textbf{Drift approximation through the control term.}
		There exists a globally Lipschitz function
		\[
		\mathbf{b}=(b_1,\dots,b_N):\mathbb{R}^{2N}\to\mathbb{R}^{2N},
		\]
		where
		\[
		b_i(x)=\theta_i(\mu_i-x_i)+\sum_{j\neq i}F_{ij}(x_i,x_j),
		\]
		with $\theta_i>0$ and $F_{ij}:\mathbb{R}^2\times\mathbb{R}^2\to\mathbb{R}^2$ globally Lipschitz, such that for some $p\ge1$,
		\[
		\mathbb{E}\!\left[
		\sup_{0\le t\le T}
		\left|
		\frac{1}{\sqrt n}\mathbf{C}^{(n)}_t
		-
		\int_0^t \mathbf{b}(\mathbf{X}^{(n)}_s)\,ds
		\right|^p
		\right]
		\longrightarrow 0.
		\]
		
		\item[(H5)] \textbf{Reflection consistency.}
		The process $\mathbf{R}^{(n)}$ is such that
		\[
		\mathbf{X}^{(n)}_t\in D^N,\qquad t\in[0,T],
		\]
		and
		\[
		\mathbf{X}^{(n)}
		=
		\Gamma\!\left(
		\mathbf{X}^{(n)}_0
		+\mathbf{M}^{(n)}
		+\frac{1}{\sqrt n}\mathbb{E}[\mathbf{A}^{(n)}-\mathbf{D}^{(n)}]
		+\frac{1}{\sqrt n}\mathbf{C}^{(n)}
		\right),
		\]
		where $\Gamma$ denotes the Skorokhod map on $D^N$.
	\end{itemize}
	
	To illustrate assumption (H2) in a concrete setting, we record the following
	componentwise functional central limit theorem for Poisson-driven input processes.
	In the vector-valued model, the convergence in (H2) is understood coordinatewise,
	and the full $2N$-dimensional noise is obtained by collecting the limits of all
	coordinates.
	
	\begin{lemma}[Componentwise FCLT for Poisson-driven input processes]
		\label{lem:fclt}
		Fix $i\in\{1,\dots,N\}$ and one coordinate of the two-dimensional process
		$Q_t^{(n),i}$. Let $A_t^{(n),i}$ and $D_t^{(n),i}$ be independent Poisson
		processes with intensities $n\lambda_i$ and $n\mu_i$, respectively. Then
		\[
		\frac{1}{\sqrt n}\left(A_t^{(n),i}-D_t^{(n),i}-n(\lambda_i-\mu_i)t\right)
		\Rightarrow
		\sigma_i W_t^i
		\quad\text{in } D([0,T],\mathbb R),
		\]
		where $W^i$ is a standard Brownian motion and $\sigma_i^2=\lambda_i+\mu_i$.
	\end{lemma}

	\begin{proof}
		Since $A_t^{(n), i}$ and $D_t^{(n), i}$ are independent Poisson processes with rates $n\lambda_i$ and $n\mu_i$, respectively, we have
		\[
		\mathbb{E}[A_t^{(n), i}] = n\lambda_i t, 
		\quad \mathrm{Var}(A_t^{(n), i}) = n\lambda_i t, 
		\qquad 
		\mathbb{E}[D_t^{(n), i}] = n\mu_i t, 
		\quad \mathrm{Var}(D_t^{(n), i}) = n\mu_i t.
		\]
		
		Define the centered and scaled processes
		\[
		M_t^{(n), i} := \frac{A_t^{(n), i} - n\lambda_i t}{\sqrt{n}}, 
		\qquad 
		N_t^{(n), i} := \frac{D_t^{(n), i} - n\mu_i t}{\sqrt{n}}.
		\]
		
		By the classical functional central limit theorem (Donsker’s theorem, see e.g., \cite{Whitt2002stochastic,Billingsley2013convergence}), as $n \to \infty$,
		\[
		M_t^{(n), i} \Rightarrow \sqrt{\lambda_i}\, W_t^{i,1}, 
		\qquad 
		N_t^{(n), i} \Rightarrow \sqrt{\mu_i}\, W_t^{i,2},
		\]
		where $W_t^{i,1}$ and $W_t^{i,2}$ are independent standard Brownian motions. 
		The convergence holds in the Skorokhod $J_1$ topology on $D([0,T], \mathbb{R})$ (see, e.g., \cite{Whitt2002stochastic,Billingsley2013convergence}).
		
		Since $(M_t^{(n), i}, N_t^{(n), i})$ converge jointly, the continuous mapping theorem \cite{Billingsley2013convergence} implies that their difference also converges:
		\[
		M_t^{(n), i} - N_t^{(n), i} 
		\Rightarrow \sqrt{\lambda_i}\, W_t^{i,1} - \sqrt{\mu_i}\, W_t^{i,2}.
		\]
		
		The limit process on the right-hand side is mean-zero Gaussian with variance
		\[
		\mathrm{Var}\!\left(\sqrt{\lambda_i}\, W_t^{i,1} - \sqrt{\mu_i}\, W_t^{i,2}\right)
		= (\lambda_i + \mu_i)t,
		\]
		since $W_t^{i,1}$ and $W_t^{i,2}$ are independent.
		
		We can therefore define
		\[
		\sigma_i W_t^i := \sqrt{\lambda_i}\, W_t^{i,1} - \sqrt{\mu_i}\, W_t^{i,2},
		\qquad \text{where } \sigma_i^2 = \lambda_i + \mu_i,
		\]
		with $W_t^i$ a standard Brownian motion. 
		Note that, since the limit process is continuous, convergence in $D([0,T],\mathbb{R})$ implies convergence in $C([0,T],\mathbb{R})$ as well.
		
		Hence,
		\[
		\frac{1}{\sqrt{n}} \left( A_t^{(n), i} - D_t^{(n), i} - n(\lambda_i - \mu_i)t \right)
		= M_t^{(n), i} - N_t^{(n), i}
		\Rightarrow \sigma_i W_t^i,
		\]
		which completes the proof.
	\end{proof}
	This lemma provides a scalar example of the primitive noise convergence assumed
	in (H2). In the present vector-valued setting, the same argument applies to each
	coordinate, and the coordinatewise limits combine to yield the $2N$-dimensional
	Brownian noise appearing in the diffusion approximation.
	
	\begin{lemma}[A priori bound for the reflected integral equation]
		\label{lem:apriori_bound}
		Let $T>0$, and let $D^N \subset \mathbb{R}^{2N}$ be a convex domain such that the associated Skorokhod map
		\[
		\Gamma : D([0,T],\mathbb{R}^{2N}) \to D([0,T],\mathbb{R}^{2N})
		\]
		is Lipschitz continuous; that is, there exists a constant $K_\Gamma>0$ such that
		\[
		\sup_{0\le t\le T} |\Gamma(\eta)(t)-\Gamma(\tilde\eta)(t)|
		\le
		K_\Gamma \sup_{0\le t\le T} |\eta(t)-\tilde\eta(t)|
		\]
		for all $\eta,\tilde\eta \in D([0,T],\mathbb{R}^{2N})$.
		
		Let $\mathbf{b}:\mathbb{R}^{2N}\to\mathbb{R}^{2N}$ be globally Lipschitz. Then there exists a constant
		$C_T>0$, depending only on $T$, $K_\Gamma$, and the Lipschitz constant of $\mathbf{b}$, such that for any càdlàg process $\mathbf{X}$ satisfying
		\begin{equation}\label{eq:reflected_IE_lemma}
			\mathbf{X}
			=
			\Gamma\!\left(
			\mathbf{x}_0+\mathbf{m}+\mathbf{a}+\int_0^\cdot \mathbf{b}(\mathbf{X}_s)\,ds+\boldsymbol{\varepsilon}
			\right)
			\quad\text{on }[0,T],
		\end{equation}
		where
		\[
		\mathbf{x}_0\in \mathbb{R}^{2N}, \qquad
		\mathbf{m},\mathbf{a},\boldsymbol{\varepsilon}\in D([0,T],\mathbb{R}^{2N}),
		\]
		one has
		\begin{equation}\label{eq:apriori_estimate}
			\sup_{0\le t\le T} |\mathbf{X}_t|
			\le
			C_T\left(
			1+|\mathbf{x}_0|
			+\sup_{0\le t\le T}|\mathbf{m}_t|
			+\sup_{0\le t\le T}|\mathbf{a}_t|
			+\sup_{0\le t\le T}|\boldsymbol{\varepsilon}_t|
			\right).
		\end{equation}
		Consequently, for every $p\ge 1$ there exists a constant $\widetilde C_{T,p}>0$ such that
		\begin{align}
			\mathbb{E}\Big[\sup_{0\le t\le T} |\mathbf{X}_t|^p\Big]
			\le
			\widetilde C_{T,p}
			\Big(
			1
			+\mathbb{E}|\mathbf{x}_0|^p
			+\mathbb{E}\Big[\sup_{0\le t\le T}|\mathbf{m}_t|^p\Big]
			+\mathbb{E}\Big[\sup_{0\le t\le T}|\mathbf{a}_t|^p\Big]
			+\mathbb{E}\Big[\sup_{0\le t\le T}|\boldsymbol{\varepsilon}_t|^p\Big]
			\Big),
			\label{eq:apriori_estimate_moment}
		\end{align}
		whenever the right-hand side is finite.
	\end{lemma}
	
	\begin{proof}
		Since $\mathbf{b}$ is globally Lipschitz, it has linear growth: there exists a constant
		$L>0$ such that
		\[
		|\mathbf{b}(x)| \le L(1+|x|),
		\qquad x\in\mathbb{R}^{2N}.
		\]
		
		Set
		\[
		\mathbf{Y}_t
		:=
		\mathbf{x}_0+\mathbf{m}_t+\mathbf{a}_t+\int_0^t \mathbf{b}(\mathbf{X}_s)\,ds+\boldsymbol{\varepsilon}_t,
		\qquad t\in[0,T].
		\]
		Then \eqref{eq:reflected_IE_lemma} may be written simply as
		\[
		\mathbf{X}=\Gamma(\mathbf{Y}).
		\]
		
		Let $\mathbf{0}$ denote the zero path in $D([0,T],\mathbb{R}^{2N})$. By the Lipschitz continuity of $\Gamma$,
		\[
		\sup_{0\le t\le T} |\mathbf{X}_t|
		=
		\sup_{0\le t\le T} |\Gamma(\mathbf{Y})(t)|
		\le
		\sup_{0\le t\le T} |\Gamma(\mathbf{Y})(t)-\Gamma(\mathbf{0})(t)|
		+\sup_{0\le t\le T}|\Gamma(\mathbf{0})(t)|.
		\]
		Hence,
		\[
		\sup_{0\le t\le T} |\mathbf{X}_t|
		\le
		K_\Gamma \sup_{0\le t\le T} |\mathbf{Y}_t|
		+\sup_{0\le t\le T}|\Gamma(\mathbf{0})(t)|.
		\]
		Absorbing the constant $\sup_{0\le t\le T}|\Gamma(\mathbf{0})(t)|$ into a generic constant $C$, we obtain
		\[
		\sup_{0\le t\le T} |\mathbf{X}_t|
		\le
		C
		+
		K_\Gamma \left(
		|\mathbf{x}_0|
		+\sup_{0\le t\le T}|\mathbf{m}_t|
		+\sup_{0\le t\le T}|\mathbf{a}_t|
		+\sup_{0\le t\le T}|\boldsymbol{\varepsilon}_t|
		+\sup_{0\le t\le T}\left|\int_0^t \mathbf{b}(\mathbf{X}_s)\,ds\right|
		\right).
		\]
		Using the linear growth of $\mathbf{b}$,
		\[
		\sup_{0\le t\le T}\left|\int_0^t \mathbf{b}(\mathbf{X}_s)\,ds\right|
		\le
		\int_0^T |\mathbf{b}(\mathbf{X}_s)|\,ds
		\le
		L\int_0^T (1+|\mathbf{X}_s|)\,ds
		\le
		LT+L\int_0^T \sup_{0\le r\le s}|\mathbf{X}_r|\,ds.
		\]
		Therefore,
		\[
		\sup_{0\le t\le T} |\mathbf{X}_t|
		\le
		C_1\left(
		1+|\mathbf{x}_0|
		+\sup_{0\le t\le T}|\mathbf{m}_t|
		+\sup_{0\le t\le T}|\mathbf{a}_t|
		+\sup_{0\le t\le T}|\boldsymbol{\varepsilon}_t|
		\right)
		+
		C_2 \int_0^T \sup_{0\le r\le s}|\mathbf{X}_r|\,ds,
		\]
		for suitable constants $C_1,C_2>0$ depending only on $T$, $K_\Gamma$, and $L$.
		
		Define
		\[
		u(t):=\sup_{0\le r\le t}|\mathbf{X}_r|, \qquad t\in[0,T].
		\]
		Then
		\[
		u(t)
		\le
		C_1\left(
		1+|\mathbf{x}_0|
		+\sup_{0\le s\le T}|\mathbf{m}_s|
		+\sup_{0\le s\le T}|\mathbf{a}_s|
		+\sup_{0\le s\le T}|\boldsymbol{\varepsilon}_s|
		\right)
		+
		C_2\int_0^t u(s)\,ds.
		\]
		By Gronwall's inequality,
		\[
		u(T)
		\le
		C_T\left(
		1+|\mathbf{x}_0|
		+\sup_{0\le t\le T}|\mathbf{m}_t|
		+\sup_{0\le t\le T}|\mathbf{a}_t|
		+\sup_{0\le t\le T}|\boldsymbol{\varepsilon}_t|
		\right),
		\]
		which proves \eqref{eq:apriori_estimate}.
		
		Finally, raising both sides of \eqref{eq:apriori_estimate} to the power $p\ge1$, using the elementary inequality
		\[
		(a_1+\cdots+a_k)^p \le C_{p,k}(a_1^p+\cdots+a_k^p),
		\]
		and then taking expectations, we obtain \eqref{eq:apriori_estimate_moment}.
	\end{proof}

	\subsubsection{Reflected SDEs and Skorokhod formulation}
	
	Let $D \subset \mathbb{R}^2$ be a convex domain with sufficiently regular boundary such that the Skorokhod problem with inward normal reflection is well posed (see, e.g., \cite{Saisho1987,Dupuis1993}).
	
	Given a continuous path $w \in C([0,T],\mathbb{R}^2)$ with $w(0)\in \bar D$, the Skorokhod problem consists in finding a pair $(\xi,\phi)$ such that
	\[
	\xi(t) = w(t) + \phi(t), \quad t\in[0,T],
	\]
	where:
	\begin{itemize}
		\item $\xi(t)\in \bar D$ for all $t\ge0$,
		\item $\phi$ is a bounded variation process with $\phi(0)=0$,
		\item $\phi$ acts only when $\xi(t)\in \partial D$ and pushes $\xi$ in the inward normal direction.
	\end{itemize}
	
	Under standard geometric conditions on $D$, this problem admits a unique solution, and the associated Skorokhod map
	\[
	\Gamma : C([0,T],\mathbb{R}^2) \to C([0,T],\mathbb{R}^2)
	\]
	is Lipschitz continuous.
	
	We now consider the reflected stochastic differential equation
	\begin{equation}\label{eq:reflected_SDE_clean}
		dX_t = f(X_t)dt + \sigma dW_t + dL_t, \qquad X_0 \in \bar D,
	\end{equation}
	where $W_t$ is a standard Brownian motion, 
	$\sigma \in \mathbb{R}^{2 \times 2}$ is a constant diffusion matrix and $L_t$ is the reflection term arising from the Skorokhod problem.
	
	\begin{definition}[Reflected SDE Solution]\label{def}
		A pair \( (X_t, L_t) \) is said to be a solution to the Skorokhod-type SDE \eqref{eq:reflected_SDE_clean}  if:
		\begin{itemize}
			\item[(i)] \( X_t \) is a continuous, \( \bar{D} \)-valued, \( \mathcal{F}_t \)-adapted process;
			\item[(ii)] \( L_t \in \mathbb{R}^2 \) is a continuous, \( \mathcal{F}_t \)-adapted process of bounded variation on every finite interval, with \( L_0 = 0 \), and
			\begin{align}
				L_t &= \int_0^t \mathbf{n}(s)\,d|L|_s, \\
				|L|_t &= \int_0^t \mathbbm{1}_{\partial D}(X_s)\, d|L|_s;
			\end{align}
			\item[(iii)] For all \( s \), the reflection direction \( \mathbf{n}(s) \in \mathcal{N}_{X_s} \), the inward normal cone to \( \partial D \) at the point \( X_s \).
		\end{itemize}
	\end{definition}

	Under Lipschitz assumptions on $f$, the reflected SDE \eqref{eq:reflected_SDE_clean} admits a unique solution (see, e.g., \cite{Saisho1987,Dupuis1993}).

	\section{Statement of the main results}
	Note that in the prelimit model, each process \(Q_t^{(n), i}\) represents the state of an individual stochastic component, while the diffusion-scaled process \(X_t^{(n), i}\) captures its centered and normalized fluctuations under diffusion scaling. The limiting dynamics therefore describe the asymptotic behavior of the interacting system in terms of a reflected stochastic differential equation.

	{\color{black}{\begin{theorem}[Diffusion approximation to an interacting reflected OU system]
				\label{thm:queue_to_sde_revised}
				Let $T>0$ and let $D^N \subset \mathbb{R}^{2N}$ be a convex domain such that the Skorokhod problem on $D^N$ with inward normal reflection is well posed and the associated Skorokhod map is Lipschitz continuous on $D([0,T],\mathbb{R}^{2N})$.
				
				For each $n\in \mathbb{N}$, consider the process
				\[
				\mathbf{Q}^{(n)}_t = \mathbf{Q}^{(n)}_0+\mathbf{A}^{(n)}_t-\mathbf{D}^{(n)}_t+\mathbf{C}^{(n)}_t+\mathbf{R}^{(n)}_t,
				\qquad t\in[0,T],
				\]
				where $\mathbf{Q}^{(n)}=(Q^{(n),1},\dots,Q^{(n),N})$ and, for each $i=1,\dots,N$,
				\[
				Q_t^{(n),i}\in \mathbb{R}^2.
				\]
				Define the diffusion-scaled process
				\[
				X_t^{(n),i}:=\frac{Q_t^{(n),i}-n\mu_i}{\sqrt n},
				\qquad
				\mathbf{X}_t^{(n)}=(X_t^{(n),1},\dots,X_t^{(n),N})\in D^N,
				\]
				where $\mu_i\in \mathbb{R}^2$.
				
				Suppose assumptions (H1)-(H5) hold.

				Then
				\[
				\mathbf{X}^{(n)} \Rightarrow \mathbf{X}
				\quad\text{in } D([0,T],\mathbb{R}^{2N}),
				\]
				where $\mathbf{X}=(X^1,\dots,X^N)$ is the unique reflected diffusion solving
				\begin{align}
					\label{eq:limit_sde_revised}
					dX_t^i
					=
					\theta_i(\mu_i-X_t^i)\,dt
					+\sum_{j\neq i}F_{ij}(X_t^i,X_t^j)\,dt
					+\Sigma_i\,dW_t^i
					+dL_t^i,
					\qquad i=1,\dots,N, 
				\end{align} with  initial condition $\mathbf{X}_0$, and where $
				\Sigma_i \in \mathbb{R}^{2\times 2}, \qquad i=1,\dots,N. 
				$,  $\mathbf{L}=(L^1,\dots,L^N)$ is the bounded-variation reflection process associated with the Skorokhod problem on $D^N$. Equivalently,
				\[
				\mathbf{X}_t
				=
				\mathbf{X}_0
				+\int_0^t \mathbf{b}(\mathbf{X}_s)\,ds
				+\Sigma \mathbf{W}_t
				+\mathbf{L}_t,
				\qquad \mathbf{X}_t\in D^N,\quad t\in[0,T].
				\]
			\end{theorem}
			
			\begin{proof}
				We divide the argument into several steps.
				
				\medskip
				\noindent
				\text{Step 1: Prelimit decomposition and representation through the Skorokhod map.}
				
				By definition,
				\[
				\mathbf{Q}^{(n)}_t
				=
				\mathbf{Q}^{(n)}_0
				+\mathbf{A}^{(n)}_t
				-\mathbf{D}^{(n)}_t
				+\mathbf{C}^{(n)}_t
				+\mathbf{R}^{(n)}_t.
				\]
				After centering at $n\boldsymbol{\mu}$ and scaling by $\sqrt n$, we obtain
				\[
				\mathbf{X}^{(n)}_t
				=
				\mathbf{X}^{(n)}_0
				+
				\frac{1}{\sqrt n}
				\Big(
				\mathbf{A}^{(n)}_t-\mathbf{D}^{(n)}_t
				\Big)
				+
				\frac{1}{\sqrt n}\mathbf{C}^{(n)}_t
				+
				\frac{1}{\sqrt n}\mathbf{R}^{(n)}_t,
				\]
				where
				\[
				\mathbf{X}^{(n)}_0
				:=
				\frac{\mathbf{Q}^{(n)}_0-n\boldsymbol{\mu}}{\sqrt n}.
				\]
				
				We further decompose the primitive noise term as
				\[
				\frac{1}{\sqrt n}
				\big(
				\mathbf{A}^{(n)}_t-\mathbf{D}^{(n)}_t
				\big)
				=
				\mathbf{M}^{(n)}_t
				+
				\frac{1}{\sqrt n}\,
				\mathbb{E}\!\left[\mathbf{A}^{(n)}_t-\mathbf{D}^{(n)}_t\right],
				\]
				where
				\[
				\mathbf{M}^{(n)}_t
				:=
				\frac{1}{\sqrt n}
				\Big(
				\mathbf{A}^{(n)}_t-\mathbf{D}^{(n)}_t
				-\mathbb{E}[\mathbf{A}^{(n)}_t-\mathbf{D}^{(n)}_t]
				\Big).
				\]
				Hence
				\begin{equation}\label{eq:prelimit-decomp}
					\mathbf{X}^{(n)}_t
					=
					\mathbf{X}^{(n)}_0
					+\mathbf{M}^{(n)}_t
					+\frac{1}{\sqrt n}\mathbb{E}\!\left[\mathbf{A}^{(n)}_t-\mathbf{D}^{(n)}_t\right]
					+\frac{1}{\sqrt n}\mathbf{C}^{(n)}_t
					+\frac{1}{\sqrt n}\mathbf{R}^{(n)}_t.
				\end{equation}
				
				By assumption (H5), the reflection mechanism is described by the Skorokhod map $\Gamma$ on $D^N$, so that
				\begin{equation}\label{eq:Sk-map-prelimit}
					\mathbf{X}^{(n)}
					=
					\Gamma(\mathbf{Y}^{(n)}),
				\end{equation}
				where
				\[
				\mathbf{Y}^{(n)}_t
				:=
				\mathbf{X}^{(n)}_0
				+\mathbf{M}^{(n)}_t
				+\frac{1}{\sqrt n}\mathbb{E}\!\left[\mathbf{A}^{(n)}_t-\mathbf{D}^{(n)}_t\right]
				+\frac{1}{\sqrt n}\mathbf{C}^{(n)}_t.
				\]
				
				\medskip
				\noindent
				\text{Step 2: Replacement of the control term by the drift integral.}
				
				By assumption (H4), for every $T>0$,
				\[
				\sup_{0\le t\le T}
				\left|
				\frac{1}{\sqrt n}\mathbf{C}^{(n)}_t
				-
				\int_0^t \mathbf{b}(\mathbf{X}^{(n)}_s)\,ds
				\right|
				\xrightarrow{\mathbb{P}} 0.
				\]
				Therefore, if we define
				\[
				\widetilde{\mathbf{Y}}^{(n)}_t
				:=
				\mathbf{X}^{(n)}_0
				+\mathbf{M}^{(n)}_t
				+\frac{1}{\sqrt n}\mathbb{E}\!\left[\mathbf{A}^{(n)}_t-\mathbf{D}^{(n)}_t\right]
				+\int_0^t \mathbf{b}(\mathbf{X}^{(n)}_s)\,ds,
				\]
				then
				\begin{equation}\label{eq:Y-tilde-close}
					\sup_{0\le t\le T}
					\left|
					\mathbf{Y}^{(n)}_t-\widetilde{\mathbf{Y}}^{(n)}_t
					\right|
					\xrightarrow{\mathbb{P}} 0.
				\end{equation}
				
				Since the Skorokhod map $\Gamma$ is Lipschitz continuous on $D([0,T],\mathbb{R}^{2N})$, it follows from \eqref{eq:Sk-map-prelimit} and \eqref{eq:Y-tilde-close} that
				\begin{equation}\label{eq:X-asymptotic-fixed-point}
					\sup_{0\le t\le T}
					\left|
					\mathbf{X}^{(n)}_t
					-
					\Gamma\!\left(
					\widetilde{\mathbf{Y}}^{(n)}
					\right)_t
					\right|
					\xrightarrow{\mathbb{P}} 0.
				\end{equation}
				
				Thus, up to an error vanishing in probability, $\mathbf{X}^{(n)}$ satisfies the reflected integral equation
				\[
				\mathbf{X}^{(n)}
				\approx
				\Gamma\!\left(
				\mathbf{X}^{(n)}_0
				+\mathbf{M}^{(n)}
				+\frac{1}{\sqrt n}\mathbb{E}[\mathbf{A}^{(n)}-\mathbf{D}^{(n)}]
				+\int_0^\cdot \mathbf{b}(\mathbf{X}^{(n)}_s)\,ds
				\right).
				\]
				
				\medskip
				\noindent
				\text{Step 3: Tightness of the sequence $\{\mathbf{X}^{(n)}\}_{n\ge1}$.}
				
				We prove that $\{\mathbf{X}^{(n)}\}$ is tight in $D([0,T],\mathbb{R}^{2N})$ by establishing uniform moment bounds and using the continuity of the Skorokhod map.
				
				\smallskip
				\noindent
				\emph{Step 3.1: Reformulation of the prelimit equation.}
				
				Define the error term
				\[
				\boldsymbol{\varepsilon}^{(n)}_t
				:=
				\frac{1}{\sqrt n}\mathbf{C}^{(n)}_t
				-
				\int_0^t \mathbf{b}(\mathbf{X}^{(n)}_s)\,ds.
				\]
				By assumption (H4),
				\[
				\mathbb{E}\Big[\sup_{0\le t\le T} |\boldsymbol{\varepsilon}^{(n)}_t|^p\Big]\to 0.
				\]
				In particular,
				\[
				\sup_{0\le t\le T} |\boldsymbol{\varepsilon}^{(n)}_t|
				\;\xrightarrow{\mathbb P}\; 0.
				\]
				
				Then the prelimit representation \eqref{eq:prelimit-decomp} can be rewritten as
				\[
				\mathbf{X}^{(n)}
				=
				\Gamma\!\left(
				\mathbf{X}^{(n)}_0
				+\mathbf{M}^{(n)}
				+\frac{1}{\sqrt n}\mathbb{E}[\mathbf{A}^{(n)}-\mathbf{D}^{(n)}]
				+\int_0^\cdot \mathbf{b}(\mathbf{X}^{(n)}_s)\,ds
				+\boldsymbol{\varepsilon}^{(n)}
				\right).
				\]
				
				\smallskip
				\noindent
				\emph{Step 3.2: Uniform moment bound.}
				
				We now apply Lemma~\ref{lem:apriori_bound} with
				\[
				\mathbf{x}_0=\mathbf{X}^{(n)}_0,\qquad
				\mathbf{m}=\mathbf{M}^{(n)},\qquad
				\mathbf{a}=\frac{1}{\sqrt n}\mathbb{E}[\mathbf{A}^{(n)}-\mathbf{D}^{(n)}],\qquad
				\boldsymbol{\varepsilon}=\boldsymbol{\varepsilon}^{(n)}.
				\]
				This yields, for some constant $C_T>0$ independent of $n$,
				\begin{align}
					\sup_{0\le t\le T} |\mathbf{X}^{(n)}_t|
					\le
					C_T\Big(
					1
					+|\mathbf{X}^{(n)}_0|
					+\sup_{0\le t\le T}|\mathbf{M}^{(n)}_t|
					+\sup_{0\le t\le T}\Big|\tfrac{1}{\sqrt n}\mathbb{E}[\mathbf{A}^{(n)}_t-\mathbf{D}^{(n)}_t]\Big|
					+\sup_{0\le t\le T}|\boldsymbol{\varepsilon}^{(n)}_t|
					\Big).
					\label{eq:uniform_bound_step3}
				\end{align}
				
				Taking $p\ge1$ and expectations, and using (H1)--(H3) together with the
				$L^p$-convergence in (H4), we obtain
				\[
				\sup_{n\ge1}
				\mathbb{E}\Big[
				\sup_{0\le t\le T} |\mathbf{X}^{(n)}_t|^p
				\Big]
				<\infty.
				\]
				
				\smallskip
				\noindent
				\emph{Step 3.3: Tightness.}
				
				By (H2), the sequence $\{\mathbf{M}^{(n)}\}$ is tight in $D([0,T],\mathbb{R}^{2N})$. By (H1), $\{\mathbf{X}^{(n)}_0\}$ is tight in $\mathbb{R}^{2N}$, and by (H3) the deterministic centering term vanishes uniformly.
				
				Furthermore, the processes
				\[
				\int_0^\cdot \mathbf{b}(\mathbf{X}^{(n)}_s)\,ds
				\]
				are tight in $C([0,T],\mathbb{R}^{2N})$ due to the uniform moment bound \eqref{eq:uniform_bound_step3} and the linear growth of $\mathbf{b}$.
				
				Finally, $\boldsymbol{\varepsilon}^{(n)}\to 0$ in probability uniformly on $[0,T]$, hence is tight.
				
				Combining these facts and using the continuity of the Skorokhod map $\Gamma$, we conclude that $\{\mathbf{X}^{(n)}\}$ is tight in $D([0,T],\mathbb{R}^{2N})$.
				
				\medskip
				\noindent
				\text{Step 4: Identification of subsequential limits.}
				
				Let $\mathbf{X}$ be any subsequential limit of $\{\mathbf{X}^{(n)}\}$. By tightness, there exists a subsequence, still denoted by $\mathbf{X}^{(n)}$, such that
				\[
				\mathbf{X}^{(n)} \Rightarrow \mathbf{X}
				\qquad\text{in } D([0,T],\mathbb{R}^{2N}).
				\]
				Together with (H1) and (H2), and after using the Skorokhod representation theorem if necessary, we may assume that on a common probability space,
				\[
				\mathbf{X}^{(n)}_0 \to \mathbf{X}_0,\qquad
				\mathbf{M}^{(n)} \to \Sigma \mathbf{W},\qquad
				\mathbf{X}^{(n)} \to \mathbf{X},
				\]
				almost surely in the Skorokhod topology. Since the limiting processes are continuous, this convergence is in fact locally uniform on $[0,T]$ almost surely.
				
				Because $\mathbf{b}$ is continuous with linear growth and $\mathbf{X}^{(n)}\to\mathbf{X}$ uniformly on $[0,T]$ almost surely, we have
				\[
				\int_0^t \mathbf{b}(\mathbf{X}^{(n)}_s)\,ds
				\to
				\int_0^t \mathbf{b}(\mathbf{X}_s)\,ds
				\qquad\text{uniformly in } t\in[0,T],
				\]
				almost surely. In addition, by (H3),
				\[
				\sup_{0\le t\le T}
				\left|
				\frac{1}{\sqrt n}\mathbb{E}\!\left[\mathbf{A}^{(n)}_t-\mathbf{D}^{(n)}_t\right]
				\right|
				\to 0.
				\]
				Therefore,
				\[
				\widetilde{\mathbf{Y}}^{(n)}
				\to
				\mathbf{Y}
				:=
				\mathbf{X}_0+\Sigma \mathbf{W}
				+\int_0^\cdot \mathbf{b}(\mathbf{X}_s)\,ds
				\]
				in $D([0,T],\mathbb{R}^{2N})$, almost surely.
				
				Applying the continuity of the Skorokhod map $\Gamma$ and using \eqref{eq:X-asymptotic-fixed-point}, we obtain
				\[
				\mathbf{X}
				=
				\Gamma(\mathbf{Y}),
				\]
				that is,
				\[
				\mathbf{X}_t
				=
				\mathbf{X}_0
				+\int_0^t \mathbf{b}(\mathbf{X}_s)\,ds
				+\Sigma \mathbf{W}_t
				+\mathbf{L}_t,
				\qquad
				\mathbf{X}_t\in D^N,
				\]
				where $\mathbf{L}$ is the reflection term associated with the Skorokhod problem on $D^N$.
				
				Equivalently, for each $i=1,\dots,N$,
				\[
				dX_t^i
				=
				\theta_i(\mu_i-X_t^i)\,dt
				+\sum_{j\neq i}F_{ij}(X_t^i,X_t^j)\,dt
				+\Sigma_i\,dW_t^i
				+dL_t^i.
				\]
				
				\medskip
				\noindent
				\text{Step 5: Uniqueness of the limit and convergence of the full sequence.}
				
				By the well-posedness assumption on the reflected SDE on $D^N$ with globally Lipschitz drift $\mathbf{b}$ and diffusion matrix $\Sigma$, the limiting reflected equation admits a unique solution in law. Hence every subsequential limit of $\{\mathbf{X}^{(n)}\}$ coincides with the same law.
				
				Therefore the entire sequence converges:
				\[
				\mathbf{X}^{(n)} \Rightarrow \mathbf{X}
				\qquad\text{in } D([0,T],\mathbb{R}^{2N}),
				\]
				where $\mathbf{X}$ is the unique reflected diffusion solving
				\[
				\mathbf{X}_t
				=
				\mathbf{X}_0
				+\int_0^t \mathbf{b}(\mathbf{X}_s)\,ds
				+\Sigma \mathbf{W}_t
				+\mathbf{L}_t.
				\]
				
				This completes the proof.
			\end{proof}
	}}

	Note that the martingale component \(\mathbf{M}_t^{(n)}\) in the semimartingale decomposition captures the intrinsic randomness of the system, arising from high-frequency stochastic fluctuations in the underlying dynamics. Under diffusion scaling, these fluctuations converge in distribution to a Brownian motion, making \(\mathbf{M}_t^{(n)}\) a key ingredient in identifying the limiting reflected diffusion. In particular, it provides the stochastic forcing in the limit and enables the application of weak convergence techniques and functional central limit theorems within the semimartingale framework.

	\paragraph{Interpretation of the diffusion scaling.}
	
	Since \(X_t^{(n), i} = \frac{Q_t^{(n), i} - n \mu_i}{\sqrt{n}}\) and \(X_t^{(n), i} \Rightarrow X_t^i\) as \(n \to \infty\), the prelimit process can be informally approximated, for large \(n\), by
	\[
	Q_t^{(n), i} \approx n \mu_i + \sqrt{n} X_t^i,
	\]
	where \(X_t^i\) is the limiting reflected interacting diffusion.
	
	This representation indicates that the process fluctuates around its nominal level \(n \mu_i\) with deviations of order \(\sqrt{n}\), consistent with classical diffusion scaling. The limiting process captures the combined effects of mean-reverting dynamics, interactions among components, and boundary constraints through reflection.
	
	From this perspective, the reflected diffusion provides a tractable continuous approximation of the underlying stochastic system, which can be useful for analyzing large-scale or constrained multi-agent dynamics motivated by queue-inspired constructions.

	\section{Numerical examples}
	
	We present numerical experiments illustrating the diffusion-scaling behavior of the stochastic systems considered in this work. Two representative scenarios are examined: (i) a crowd dynamics setting and (ii) a neural population model. In both cases, we compare discrete stochastic systems with their corresponding reflected diffusion models. The goal is not to verify convergence in a strict theoretical sense, but rather to demonstrate that, as the system size increases, the discrete dynamics exhibit behavior consistent with the diffusion approximation.
	
	\subsection{Human crowd dynamics}
	
	We consider a system of interacting agents evolving in the two-dimensional square domain
	\[
	\mathcal{D} = [0,2.5] \times [0,2.5].
	\]
	All agents share the same target location (exit), located at
	\[
	\mu = (1.25, 0).
	\]
	Three agents are initialized at
	\[
	x_1(0) = (0.75, 0.75), \quad
	x_2(0) = (1.75, 1.75), \quad
	x_3(0) = (1.25, 0.25).
	\]
	
	Each agent evolves according to a reflected stochastic dynamics consisting of:
	(i) a mean-reverting drift $\theta_i(\mu - x_i)$ with $\theta_i = 3.0$,
	(ii) additive Brownian noise with intensity $\sigma_i = 0.5$, and
	(iii) pairwise repulsive interactions modeled through a softened inverse-square force to prevent overlap.
	
	{\color{black}{
			Reflection at the boundary of $\mathcal D$ is implemented numerically through a penalty-type confinement mechanism that prevents agents from leaving the admissible domain. This approach is adopted for computational simplicity and robustness and should be viewed as a numerical approximation of the reflected interacting diffusion considered in Theorem~\ref{thm:queue_to_sde_revised}. Consequently, the simulations are intended to illustrate the diffusion-scaling behavior predicted by the theoretical framework rather than to provide a direct numerical verification of the convergence theorem. In particular, the focus is on the emergence of mean-reverting motion, interaction effects, and confinement as the system size increases.
	}}

	\medskip
	\noindent
	\textbf{Discrete stochastic approximation.}
	To approximate the diffusion dynamics, we simulate a family of $n$-dependent discrete stochastic systems whose scaled states are given by
	\[
	X_i^{(n)}(t) = \frac{Q_i^{(n)}(t) - n \mu_i}{\sqrt{n}},
	\]
	where $Q_i^{(n)}$ evolves through state-dependent stochastic increments designed to mimic the drift, interaction, and confinement effects of the diffusion model. {\color{black}{This construction can be interpreted as a diffusion-scaled birth-death-type
			stochastic approximation and is intended only as a representative prelimit
			stochastic system rather than a classical queueing network.}}

	We consider system sizes $n \in \{50, 200, 800\}$ and perform $1000$ independent simulation runs for each value of $n$. Empirical averages of the scaled processes are compared with trajectories of the corresponding diffusion model.
	
	\begin{figure}[h!]
		\centering
		\includegraphics[width=0.95\textwidth]{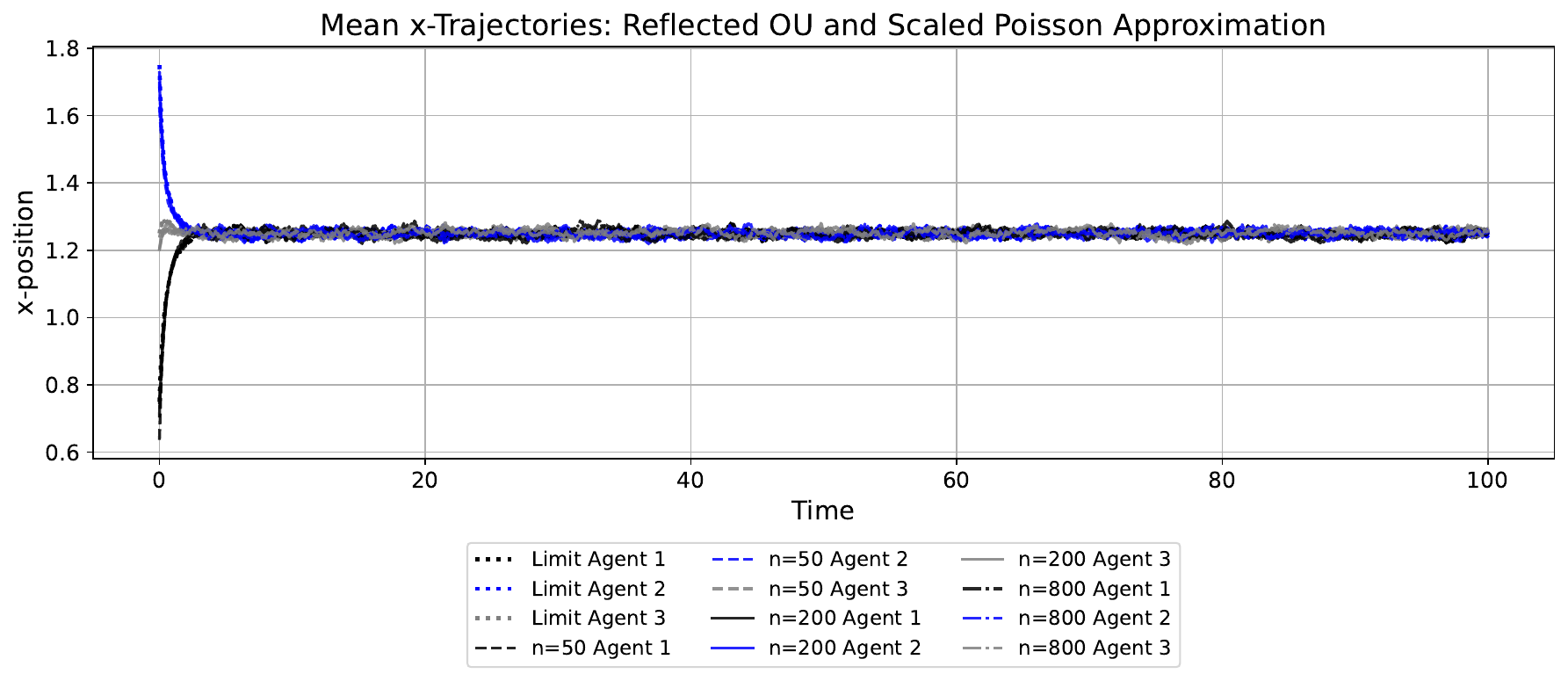}
		\caption{[Color online] Crowd dynamics under diffusion scaling.
			Mean $x$-coordinate trajectories of three interacting agents obtained from
			the reflected diffusion model (dotted lines) and the diffusion-scaled
			stochastic approximations for system sizes $n=50,200,800$. As the scaling
			parameter increases, the empirical mean trajectories of the stochastic
			systems exhibit improved agreement with the corresponding reflected
			diffusion trajectories, illustrating the diffusion-scaling behavior
			predicted by the theoretical framework.}
		\label{fig:1}
	\end{figure}
	
	Figure~\ref{fig:1} shows the mean trajectories of the $x$-coordinates of the agents. As $n$ increases, the empirical averages of the discrete systems exhibit
	improved agreement with the trajectories of the diffusion model, while
	smaller values of $n$ display more pronounced stochastic fluctuations.
	
	%
	
	\begin{table}[h!]
		\centering
		\caption{Empirical mean trajectory discrepancy between the scaled jump
			process and the diffusion model in the crowd dynamics example.}
		\label{tab:mse-crowd}
		\begin{tabular}{cc}
			\hline
			System size $n$ & Mean trajectory MSE \\
			\hline
			$50$  & $1.9192\times 10^{-4}$ \\
			$200$ & $9.9000\times 10^{-5}$ \\
			$800$ & $9.5090\times 10^{-5}$ \\
			\hline
		\end{tabular}
	\end{table}
	
	{\color{black}{The empirical discrepancy between the mean scaled jump trajectories and
			the mean diffusion trajectories is reported in Table~\ref{tab:mse-crowd}.
			The mean trajectory MSE decreases from $1.9192\times 10^{-4}$ for
			$n=50$ to $9.509\times 10^{-5}$ for $n=800$, illustrating improved
			agreement under diffusion scaling. We emphasize that this table is used
			only as a qualitative numerical illustration, since the theoretical result
			establishes weak convergence in distribution rather than strong pathwise
			or mean-square convergence.}}
	
	Overall, this experiment demonstrates that the discrete stochastic systems capture the qualitative behavior of the reflected diffusion model, including mean-reverting motion, interaction effects, and confinement within the domain.
	
	\subsection{Diffusion scaling in neural population dynamics}
	
	We next consider a system of five interacting neural populations evolving in a two-dimensional activity space. Each population is represented by a state variable $R_t^i \in \mathbb{R}^2$, with initial conditions
	\[
	R_0^1 = (1.0, 1.0), \quad
	R_0^2 = (2.0, 4.0), \quad
	R_0^3 = (3.5, 1.5), \quad
	R_0^4 = (4.5, 3.0), \quad
	R_0^5 = (2.5, 2.5).
	\]
	
	The populations evolve under:
	\begin{itemize}
		\item[(i)] a drift toward a common attractor $\mu_i = (2.5, 0.0)$,
		\item[(ii)] mutual inhibitory interactions represented by functions $F_{ij}$, and
		\item[(iii)] reflective constraints that confine the dynamics within the domain $[0,5]^2$.
	\end{itemize}
	
	{\color{black}{
			As in the crowd dynamics example, confinement within the admissible activity domain is implemented numerically through a penalty-type mechanism rather than the exact Skorokhod reflection used in the theoretical analysis. Therefore, the numerical experiments should be interpreted as illustrative examples of the diffusion-scaling regime described by Theorem~\ref{thm:queue_to_sde_revised}. Their purpose is to demonstrate qualitatively how the behavior of the finite stochastic systems approaches that of the corresponding diffusion model as the system size increases, rather than to provide a rigorous numerical validation of the convergence result.
	}}

	\medskip
	\noindent
	\textbf{Finite-size stochastic system.}
	We construct a discrete stochastic approximation in which each population evolves through random increments depending on its current state and interactions with other populations. Under diffusion scaling, the centered and normalized processes
	\[
	X_t^{(n),i} = \frac{Q_t^{(n),i} - n \mu_i}{\sqrt{n}}
	\]
	are expected to approximate the diffusion system as $n$ becomes large.
	
	We simulate the system for $n = 50, 200, 800$, using Monte Carlo sampling to estimate mean trajectories and discrepancies relative to the diffusion model.
	
	\begin{figure}[h!]
		\centering
		\includegraphics[width=0.95\textwidth]{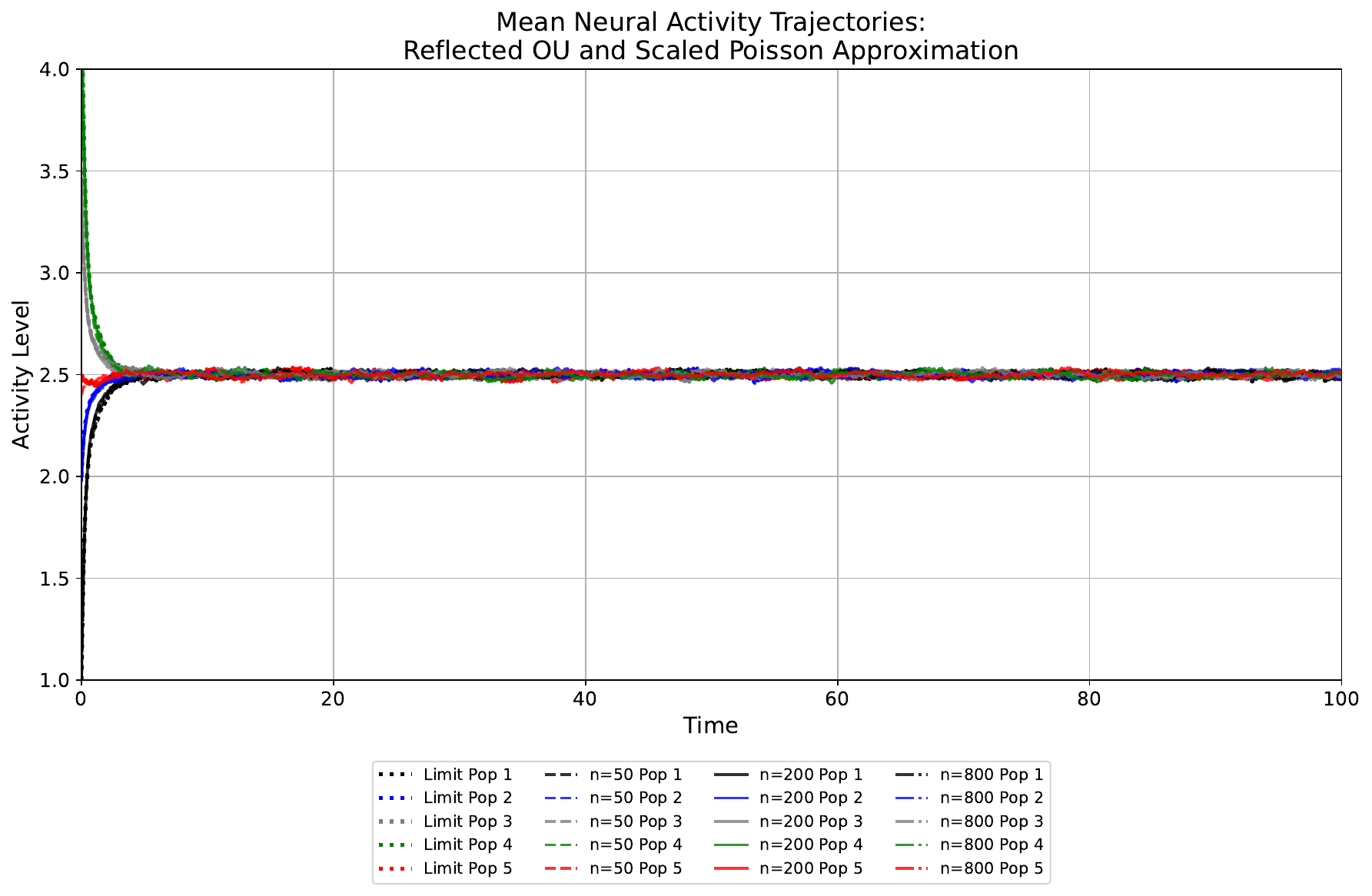}
		\caption{[Color online] Neural population dynamics under diffusion scaling.
			Mean activity trajectories of five interacting neural populations obtained
			from the reflected diffusion model (dotted lines) and the diffusion-scaled
			stochastic approximations for system sizes $n=50,200,800$. The finite
			stochastic systems display increasingly similar behavior to the reflected
			diffusion model as the scaling parameter increases, illustrating the
			emergence of diffusion-scale dynamics in interacting neural populations.}
		\label{fig:2}
	\end{figure}
	
	Figure~\ref{fig:2} illustrates the mean trajectories of the populations. As the system size increases, the discrete dynamics exhibit reduced variability and increasingly resemble the diffusion trajectories.

	\begin{table}[h!]
		\centering
		\caption{Empirical mean trajectory discrepancy between the scaled neural
			jump process and the diffusion model.}
		\label{tab:mse-neural}
		\begin{tabular}{cc}
			\hline
			System size $n$ & Mean trajectory MSE \\
			\hline
			$50$  & $1.8802\times 10^{-4}$ \\
			$200$ & $1.5242\times 10^{-4}$ \\
			$800$ & $1.4351\times 10^{-4}$ \\
			\hline
		\end{tabular}
	\end{table}
	
	{\color{black}{Table~\ref{tab:mse-neural} reports the empirical discrepancy between the
			mean scaled neural jump trajectories and the mean diffusion trajectories.
			The mean trajectory MSE decreases as $n$ increases, illustrating improved
			agreement under diffusion scaling. }}
	
	This example demonstrates that reflected diffusion models provide a useful approximation framework for describing the collective behavior of interacting stochastic populations, even in the presence of nonlinear interactions and boundary constraints.
	
	%
	
	\section{Conclusions and discussion}
	
	In this paper, we developed a diffusion-approximation framework for a class of interacting stochastic systems with reflection and control. Starting from discrete stochastic processes that incorporate stochastic fluctuations, interaction effects, and boundary constraints, we established conditions under which the rescaled dynamics converge in distribution to systems of interacting reflected Ornstein-Uhlenbeck processes. The limiting model captures key features of the underlying systems, including mean-reverting behavior, inter-agent interactions, and confinement within bounded domains.
	
	{\color{black}{The framework is motivated by diffusion-scaling limits of interacting stochastic systems with reflection and control. Unlike mean-field and McKean-Vlasov approaches, which focus on asymptotic limits as the number of interacting particles tends to infinity, the present work studies diffusion approximations arising from the scaling of underlying stochastic dynamics and leads to finite-dimensional systems of interacting reflected diffusions.
			
			We illustrated the approach through two representative numerical examples: a crowd dynamics scenario and a neural population model. The simulations demonstrate qualitative agreement between the finite stochastic systems and the corresponding reflected diffusion models. In addition, empirical discrepancy measurements show improved agreement under diffusion scaling as the scaling parameter increases. These examples illustrate how diffusion approximations can capture important qualitative features of interacting stochastic systems, including goal-directed motion, interaction effects, and boundary constraints.}}
	
	In the context of neural population dynamics, the results suggest that stochastic systems with discrete event-driven behavior may be approximated by continuous reflected stochastic differential equations at an appropriate scale. Such approximations provide a useful framework for studying population-level activity in the presence of interactions and confinement effects. Potential applications include models of competitive decision-making, sensorimotor representations, and constrained cortical dynamics.
	
	Several directions for future research remain. These include extending the analysis to more general interaction structures, time-dependent or irregular domains, and heterogeneous agent dynamics. Another important direction is the incorporation of control and optimization mechanisms within the reflected diffusion framework. It would also be of interest to investigate diffusion approximations for systems with mean-field interactions and to establish analogous results in more general reflected stochastic settings.

  \bibliographystyle{elsarticle-num} 

\bibliography{mybibn}

\begin{thebibliography}{10}
\expandafter\ifx\csname url\endcsname\relax
  \def\url#1{\texttt{#1}}\fi
\expandafter\ifx\csname urlprefix\endcsname\relax\def\urlprefix{URL }\fi
\expandafter\ifx\csname href\endcsname\relax
  \def\href#1#2{#2} \def\path#1{#1}\fi

\bibitem{Skorohod1961}
A.~V. Skorohod, Stochastic equations for diffusion process in a bounded domain,
  Theory of Probability and Its Applications VI (1961) 264--274.

\bibitem{Lions1984}
P.~L. Lions, Stochastic differential equations with reflecting boundary
  conditions, Communications on Pure and Applied Mathematics XXXVII (1984)
  511--537.

\bibitem{Saisho1987}
Y.~Saisho, Stochastic differential equations for multi-dimensional domain with
  reflecting boundary, Probab. Th. Rel. Fields 74 (1987) 455--477.

\bibitem{Dupuis1993}
P.~Dupuis, H.~Ishii, Sdes with oblique reflection on nonsmooth domains,
  Mathematical and Computer Modeling 1 (1993) 554--580.

\bibitem{Slominski1993}
L.~Slominski, On existence, uniqueness and stability of solutions of
  multidimensional {SDE}’s with reflecting boundary conditions, Ann. Inst.
  Henri. Poincar\'{e} 29 (1993) 163--198.

\bibitem{Slominski1994}
L.~S\l{l}omi\'{n}ski, On approximation of solutions of multidimensional {SDE's}
  with reflecting boundary conditions, Stochastic Processes and their
  Applications 50~(2) (1994) 197--219.

\bibitem{MarinRubio2004}
P.~Mar\'{i}n-Rubio, J.~Real, Some results on stochastic differential equations
  with reflecting boundary conditions, Journal of Theoretical Probability 17
  (2004) 705--716.

\bibitem{Wells2006}
C.~G. Wells, A stochastic approximation scheme and convergence theorem for
  particle interactions with perfectly reflecting boundary conditions, Monte
  Carlo Methods and Applications 12~(3) (2006) 291--342.

\bibitem{Sabelfeld2019}
K.~Sabelfeld, Stochastic algorithm for solving transient diffusion equations
  with a precise accounting of reflection boundary conditions on a substrate
  surface, Appl. Math. Letters 96 (2019) 187--194.

\bibitem{Li2020novel}
Y.~Li, H.~Chen, M.~Feng, A novel model for the traffic of urban roads based on
  queuing theory, in: 2020 International Conference on Intelligent Computing,
  Automation and Systems (ICICAS), IEEE, 2020, pp. 190--194.

\bibitem{Bhattacharya2021random}
R.~N. Bhattacharya, E.~C. Waymire, Random Walk, Brownian Motion, and
  Martingales, Vol.~52, Springer, 2021.

\bibitem{Falkowski2025}
A.~Falkowski, {SDE}s with two reflecting barriers driven by optional processes
  with regulated trajectories, Stochastic Processes and their Applications 179
  (2025) 104509.

\bibitem{Bass2025uniqueness}
R.~F. Bass, K.~Burdzy, Uniqueness for the {S}korokhod problem in an orthant:
  {C}ritical cases, Electronic Journal of Probability 30 (2025) 1--19.

\bibitem{Carmona2013control}
R.~Carmona, F.~Delarue, A.~Lachapelle, Control of {McKean--Vlasov} dynamics
  versus mean field games, Mathematics and Financial Economics 7~(2) (2013)
  131--166.

\bibitem{Liu2021long}
W.~Liu, L.~Wu, C.~Zhang, Long-time behaviors of mean-field interacting particle
  systems related to {McKean--Vlasov} equations, Communications in Mathematical
  Physics 387~(1) (2021) 179--214.

\bibitem{Jiang2026learning}
Q.~Jiang, L.~Li, L.~Zhang, L.~Wan, Learning collective multicellular dynamics
  with an interacting mean field neural {SDE} model, PLOS Computational Biology
  22~(1) (2026) e1013916.

\bibitem{Bramson2006stability}
M.~Bramson, Stability and heavy traffic limits for queueing networks,
  Proceedings of the XXXVIth International Probability Summer School,
  Saint-Flour, France. Springer-Verlag (2006).

\bibitem{Lee2011}
C.~Lee, A.~Weerasinghe, Convergence of a queueing system in heavy traffic with
  general patience-time distributions, Stochastic Processes and their
  Applications 121~(11) (2011) 2507--2552.

\bibitem{Pilipenko2014introduction}
A.~Pilipenko, An Introduction to Stochastic Differential Equations with
  Reflection, Universit{\"a}tsverlag Potsdam, 2014.

\bibitem{Kumaran2019queuing}
S.~K. Kumaran, D.~P. Dogra, P.~P. Roy, Queuing theory guided intelligent
  traffic scheduling through video analysis using dirichlet process mixture
  model, Expert systems with applications 118 (2019) 169--181.

\bibitem{Ata2024singular}
B.~Ata, J.~M. Harrison, N.~Si, Singular control of (reflected) brownian motion:
  a computational method suitable for queueing applications, Queueing Systems
  108~(3) (2024) 215--251.

\bibitem{Xu2024pigat}
Q.~Xu, Y.~Pang, X.~Zhou, Y.~Liu, Pigat: Physics-informed graph attention
  transformer for air traffic state prediction, IEEE Transactions on
  Intelligent Transportation Systems 25~(9) (2024) 12561--12577.

\bibitem{Su2025improved}
W.~Su, C.~Mu, L.~Xue, X.~Yang, S.~Zhu, An improved traffic coordination control
  integrating traffic flow prediction and optimization, Engineering
  Applications of Artificial Intelligence 143 (2025) 109969.

\bibitem{Vandaele2000queueing}
N.~Vandaele, T.~Van~Woensel, A.~Verbruggen, A queueing based traffic flow
  model, Transportation Research Part D: Transport and Environment 5~(2) (2000)
  121--135.

\bibitem{Helbing2003section}
D.~Helbing, A section-based queueing-theoretical traffic model for congestion
  and travel time analysis in networks, Journal of Physics A: Mathematical and
  General 36~(46) (2003) L593.

\bibitem{Bogachev2024approximate}
M.~I. Bogachev, N.~S. Pyko, N.~Tymchenko, S.~A. Pyko, O.~A. Markelov,
  Approximate waiting times for queuing systems with variable cross-correlated
  arrival rates, Physica A: Statistical Mechanics and its Applications 654
  (2024) 130152.

\bibitem{Whitt2002stochastic}
W.~Whitt, Stochastic-process Limits: An Introduction to Stochastic-process
  Limits and Their Application to Queues, Springer, 2002.

\bibitem{Billingsley2013convergence}
P.~Billingsley, Convergence of Probability Measures, John Wiley \& Sons, 2013.

\end{thebibliography}



%
%
%
\end{document}